\documentclass[12pt]{article}

\usepackage{amssymb}
\usepackage{amsmath,amsthm}
\usepackage{latexsym}
\usepackage{amsbsy}
\usepackage{amscd}
\usepackage{color}

\usepackage[all]{xy}
\CompileMatrices

\newtheorem{teorema}{Theorem}[section]
\newtheorem{proposicion}[teorema]{Proposition}
\newtheorem{lema}[teorema]{Lemma}
\newtheorem{definicion}[teorema]{Definition}
\newtheorem{corolario}[teorema]{Corollary}
\newtheorem{nota}[teorema]{Remark}

\newtheorem{ejemplo}[teorema]{Example}
\newtheorem{ejemplos}[teorema]{Examples}

\def\cajita{$\framebox{\begin{picture}(.7,.7)\put(0,0){}\end{picture}}$}
\def\gen{\mathfrak{genus}}

\def\ev{\mathrm{ev}}
\def\immto{\looparrowright}
\def\P{\mathrm{P}}
\def\rojo{} %\textcolor{red}}
\def\colorado{} %\textcolor{red}}

\title{Symmetric topological complexity \\ of projective and lens spaces}
\author{Jes\'us Gonz\'alez and Peter Landweber}
\date{\empty}

\begin{document}
\maketitle

\vspace{-3mm}
\centerline{\it Dedicated to the memory of Bob Stong}
\medskip 

\begin{abstract}
For real projective spaces, (a) the 
Euclidean immersion dimension, (b) the existence of axial maps, and
(c) the topological complexity
are known to be three facets of the same problem. 
\rojo{But when it comes to embedding dimension, the classical 
work of Berrick, Feder and Gitler leaves a
small indeterminacy when trying to identify the existence
of Euclidean embeddings of these manifolds with the existence 
of symmetric axial maps. As an alternative we show that the 
symmetrized version of (c) captures, in a sharp way,
the embedding problem}. Extensions to the case of \rojo{even}-torsion
lens spaces and complex projective spaces are  discussed.
\end{abstract}



\section{Main result}\label{mainsec}
The Euclidean immersion and embedding questions for projective 
spaces were topics of intense research during the beginning of 
the second half of the last century. In the case of real projective 
spaces, the immersion problem has recently received a fresh push, 
partly in view of a surprising reformulation in terms of a basic 
concept arising in robotics, namely, the motion planning problem 
of mechanical systems. In more detail, it is shown in~\cite{FYT} that for
$r\neq 1,3,7$, the immersion dimension of $\P^r$---the $r$-dimensional 
real projective space---agrees with $\mbox{TC}(\P^r)-1$, one unit less 
than the topological complexity  of $\P^r$
(see Definition~\ref{TCdef} and Theorem~\ref{recuperado} 
below). In this paper we accomplish a 
completely analogous goal by connecting the Euclidean embedding 
dimension of $\P^r$ with Farber-Grant's notion of symmetric motion planning. 
Before stating our main results, we recall the relevant definitions.

\medskip
The Schwarz genus (\cite{schwarz}) of a fibration $p\colon E\to B$, 
denoted by $\gen(p)$, 
is the smallest number of open
sets $U$ covering $B$ in such a way that $p$ 
admits a (continuous) section over each $U$. 

\begin{definicion}\label{TCdef}{\em
The topological complexity of 
a space $X$, TC$(X)$, is defined as the genus of the end-points evaluation
map $\ev\colon P(X)\to X\times X$, where $P(X)$ is the free path
space $X^{[0,1]}$ with the compact-open topology. 
}\end{definicion}

TC$(X)$ is a homotopy invariant
of $X$. Thinking of $X$ as the space of configurations of a given 
mechanical system, TC$(X)$ gives a measure of the topological instabilities 
in a motion planning algorithm for $X$---a perhaps discontinuous (but global) 
section of the map $\ev$. We refer the reader to \cite{Farbersur} for
a very useful survey of results in this area, \rojo{and to the 
book~\cite{farberinvitation} for 
a thorough introduction to the new mathematical
discipline of Topological Robotics.}

\medskip
We now come to the main definition (introduced and explored in~\cite{FGsymm}).
For a topological space $X$, let $\Delta_X$ be the diagonal in $X\times X$, 
and $\,\ev_1\colon P_1(X)\to X\times X-\Delta_X$ be the restriction of 
the fibration $\ev$
in Definition~\ref{TCdef}. Thus $P_1(X)$ is the subspace of $P(X)$ 
consisting of paths $\gamma : [0,1] \to X$ with $\gamma (0) \neq \gamma (1)$.
Note that $\ev_1$ is a $\mathbb{Z}/2$-equivariant map, 
where $\mathbb{Z}/2$
acts freely on both $P_1(X)$ and $X\times X-\Delta_X$, by 
running a path backwards in the former, and by 
switching coordinates in the latter. Let $P_2(X)$
and $B(X,2)$ denote the corresponding orbit spaces,
and let $\ev_2\colon P_2(X)\rightarrow B(X,2)$ denote
the resulting fibration. 

\begin{definicion}\label{TCSdef}{\em
With the above conditions, the 
symmetric topological complexity of $X$, $\mbox{TC}^S(X)$, is defined by
$\mbox{TC}^S(X)=\gen(\ev_2)+1$.
}\end{definicion}

Finally, let $E(r)$ stand for the Euclidean embedding dimension of $\P^r$.
Then, our main result is:

\begin{teorema}\label{mainT}
For $r>15$ as well as for $r\in\{1,2,4,8,9,13\}$, the symmetric topological 
complexity of $\P^r$ satisfies $\mbox{\em TC}^S(\P^r)=E(r)+1$.
\end{teorema}

Before introducing an alternative characterization of $\mbox{TC}^S(\P^r)$
(which implies Theorem~\ref{mainT}), it is convenient to say
a few words (to be expanded in Section~\ref{relax}) 
comparing the immersion and the embedding situations. 
It is known that the key concept bridging the 
immersion dimension of $\P^r$ to its topological complexity is 
that of an axial map. Not only do axial maps capture, in a 
sharp way, the immersion problem for projective 
spaces~(\cite{AGJ:axial}), but as shown in~\cite{FYT}, they 
conveniently encode instructions for the motion planning problem associated 
to $\mbox{TC}(\P^r)$. 
Now, the work in~\cite{BFG} does show a relation, at least in Haefliger's 
metastable range, between embeddings of real projective spaces on the one 
hand, and symmetric axial maps on the other. However, as of today, this 
relation has an unsettled indeterminacy of one dimension---spelled out 
in~(\ref{Eas}) below. Instead, motivated by the main trick 
in~\cite{BFG} (see the proof of Proposition~\ref{trick}), our approach
leads to a direct proof of Theorem~\ref{mainT}. To this 
end we actually need to give up using the concept of symmetric axial map, 
and replace it by that of the level of an involution (as defined 
in~(\ref{leveldef}) below). This allows us to get the following sharp and
unrestricted characterization for the symmetric topological complexity of 
$\P^r$. 

\begin{teorema}\label{Tfinal}
\rojo{For all values of $r$, $\mbox{\em TC}^S(\P^r)=
\mbox{\em level}
\left(\P^r\times\P^r-\Delta_{\P^r},\mathbb{Z}/2\right)+1$.}
\end{teorema}

Here the pair $\left(\P^r\times\P^r-\Delta_{\P^r},\mathbb{Z}/2\right)$
stands for the $\mathbb{Z}/2$-action on $\P^r\times\P^r-\Delta_{\P^r}$ that 
interchanges coordinates, whereas the level of a principal 
$\mathbb{Z}/2$-action on a space $X$, $\mbox{level\hspace{.4mm}}
(X,\mathbb{Z}/2)$, is defined by the formula
\begin{equation}\label{leveldef}
\mbox{level\hspace{.4mm}}(X,\mathbb{Z}/2)=\min\{\ell>0\,\colon\exists
\mbox{ $\mathbb{Z}/2$-equivariant map } X\to S^{\ell-1}\}
\end{equation}
where the sphere is considered with the antipodal $\mathbb{Z}/2$-action
(see~\cite{dailam}).

\medskip
The paper is organized as follows. Section~\ref{mainproof} is devoted to the 
proof of Theorem~\ref{Tfinal}. After observing that Theorem~\ref{mainT} is a 
consequence of Theorem~\ref{Tfinal} together with 
Haefliger's characterization of 
Euclidean embeddings of smooth manifolds, in Section~\ref{metastable} we make 
an {\it ad hoc} analysis of the numerical values of $\mbox{TC}^S
(\P^r)$ for those cases of $r$ outside Haefliger's metastable range.
Section~\ref{relax} surveys the relation of axial maps to immersion 
dimension~(\cite{AGJ:axial}), and to (non-necessarily symmetric)
topological complexity~(\cite{FYT}), focusing on 
the way those ideas compare to (and motivate)
our results. In Section~\ref{lens} we study the 
symmetric topological complexity of $m$-torsion lens spaces,
for $m$ even. Here our results are weaker than the case 
\rojo{$m=2$}, due in part to the fact that, as the \rojo{$2$}-torsion 
increases, the end terms in~(\ref{sinsym}) below
start measuring different phenomena, thus preventing us from closing the 
cycle of inequalities. Yet, we manage to give alternative 
characterizations (Theorem~\ref{caracte} and Proposition~\ref{symigualdade})
for the symmetric topological complexity of even-torsion 
lens spaces. One of these characterizations leads to
a particularly convenient upper bound
(Corollary~\ref{cota}) which depends not only on the dimension of the lens 
space, but also on its torsion. Theorem~\ref{stable1} 
illustrates the use of such an upper
bound in the context of non-symmetric topological complexity. 
In the final Section~\ref{cx} we compute
the numerical value of the symmetric topological complexity of 
complex projective spaces.

\medskip
\rojo{The first author gratefully acknowledges the kind support received 
from Professor
Michael Farber and the Department of Mathematical Sciences at 
Durham University during a visit in October 2008. Michael Farber's suggestions
to an earlier version of this work were very helpful. In particular, 
Farber noticed that our proof of Theorem~\ref{mainT} 
actually leads to the proof of Theorem~\ref{Tfinal}, and that 
our results and methods in Section~\ref{lens} 
(originally written for $2^e$-torsion) 
apply just as well to even-torsion lens spaces. Helpful suggestions of an 
anonymous referee to an earlier version of this work lead to a substantial
improvement in the organization of the paper.}

\section{Main proof}\label{mainproof}
There are three ingredients in the proof of Theorem~\ref{Tfinal}. For the first 
one we note that Corollary~1 on page~97 of~\cite{schwarz} 
affirms that the canonical 
$\mathbb{Z}/2$-cover $S^{n-1}\to \P^{n-1}$ classifies 
$\mathbb{Z}/2$-covers of genus at most $n$. Explicitly, the principal 
$\mathbb{Z}/2$-actions on a space $X$ which admit a 
$\mathbb{Z}/2$-equivariant map $X\to S^{n-1}$ are precisely those for which
the canonical projection $\rojo{p}\colon\!X\to X/(\mathbb{Z}/2)$ has 
$\gen\le n$. In particular, 
\begin{equation}\label{class}
\gen(p)=\mbox{level\hspace{.4mm}}(X,\mathbb{Z}/2).
\end{equation}
Propositions~\ref{symigualdad} and~\ref{trick} below
are the other two auxiliary ingredients. They are based on
the following preliminary constructions. For a path $\gamma\in P(\P^r)$, let 
$\widehat{\gamma}\colon 
[0,1]\to S^r$ be any lifting of $\gamma$ through the canonical projection 
$S^r\to \P^r$, and then set $f(\gamma)$ to be the class of 
$(\widehat\gamma(0),\widehat\gamma(1))$ in the Borel construction 
$S^r\times_{\mathbb{Z}/2}S^r=\left(S^r\times S^r\right)/ (-x,y)\sim(x,-y)$. This
gives a $\mathbb{Z}/2$-equivariant commutative diagram
\begin{equation}\label{triangulo}
\raisebox{1.5mm}{
\begin{picture}(0,22)(65,-17)
\put(0,0){
\xymatrix{
{ P(\P^r) } \ar[rr]^f \ar[rd]_{\ev} & 
& {S^r\times_{\mathbb{Z}/2} S^r } \ar[ld]^{\pi} \\
 & { \P^r\times \P^r }  & 
}}\end{picture}
}\end{equation}

\vspace{4mm}\noindent
where $\pi$ is the canonical projection, and 
the $\mathbb{Z}/2$-action on $S^r\times_{\mathbb{Z}/2} S^r$ switches 
coordinates (and the $\mathbb{Z}/2$-actions on $P(\P^r)$ and $\P^r\times \P^r$
are the obvious extensions of the respective $\mathbb{Z}/2$-actions 
on $P_1(\P^r)$ and $\P^r\times \P^r-\Delta_{\P^r}$ described  
just before Definition~\ref{TCSdef}). In particular, 
by restricting to $\P^r\times \P^r-\Delta_{\P^r}$ and then passing to 
$\mathbb{Z}/2$-orbit spaces,~(\ref{triangulo}) yields corresponding triangles

\begin{equation}\label{triangulo1}
\raisebox{2mm}{
\begin{picture}(0,21)(137,-19)
\put(0,0){
\xymatrix{
{ P_1(\P^r) } \ar[rr]^{f_1} \ar[rd]_{\ev_1} & 
& { E_1 } \ar[ld]^{\pi_1} \\
 & { \P^r\times \P^r-\Delta_{\P^r} }  & 
}\qquad
\xymatrix{
{ P_2(\P^r) } \ar[rr]^{f_2} \ar[rd]_{\ev_2} & 
& { E_2 } \ar[ld]^{\pi_2} \\
 & { B(\P^r,2) }  & 
}}\end{picture}
}\end{equation}

\vspace{2.7mm}
\begin{proposicion}\label{symigualdad}
For $i\in\{1,2\}$, $\gen(\ev_i)=\gen(\pi_i)$.
\end{proposicion}

\begin{proof}It 
suffices to construct a fiber-preserving $\mathbb{Z}/2$-equivariant 
map $g_1\colon E_1\to {P_1}(\P^r)$ running backwards in the
left triangle of~(\ref{triangulo1}).
To this end, we use a straightforward adaptation of the idea in the
first part of the proof of~\cite[Proposition~17]{FYT}. 
An explicit model for $E_1$ is the quotient of $S^r\times
S^r-{\widetilde{\Delta}}$ by the relation $(x,y)\sim(-x,-y)$, where
${\widetilde\Delta}\subset S^r \times S^r$ is given by
$\widetilde\Delta=\{(x,y)\in S^r\times S^r \,|\, x\neq\pm y\}$.
In these terms, the
$\mathbb{Z}/2$-action on $E_1$ interchanges coordinates. 
Then, the required map $g_1$ takes the class of a pair $(x_1,x_2)$
into the curve $[0,1]\to S^r\to \P^r$, where the second map 
is the canonical projection,
and the first map is given by 
\begin{equation}\label{curva}
t\mapsto \nu(tx_2+(1-t)x_1). 
\end{equation}
Here $\nu\colon \mathbb{R}^{r+1}-\{0\}\to S^r$ is the normalization map.
\end{proof}

The main trick in~\cite{BFG} is adapted for the proof of the following result.

\begin{proposicion}\label{trick}
If $\rho\colon \P^r\times\P^r-\Delta_{P^r}\to B(\P^r,2)$ stands for 
the canonical projection associated to the involution $(\P^r\times\P^r-
\Delta_{P^r},\mathbb{Z}/2)$ in Theorem~\ref{Tfinal}, then $\gen(\rho)=
\gen(\pi_2)$.
\end{proposicion}

\begin{proof}
As indicated in the proof of Proposition~\ref{symigualdad},
$E_2$ is the quotient of $S^r\times S^r-\widetilde\Delta$ by the relations
\begin{equation}\label{reltcs}
(-x,-y)\sim(x,y)\sim(y,x).
\end{equation}
Likewise, the space 
$\P^r\times\P^r-\Delta_{\P^r}$ is the quotient of $S^r\times S^r-\widetilde
\Delta$ by the relations
\begin{equation}\label{relemb}
(-x,y)\sim(x,y)\sim(x,-y).
\end{equation}
Moreover, the map
\begin{equation}\label{extendida}
S^r\times S^r-\widetilde\Delta\stackrel{\Psi}{\longrightarrow} 
S^r\times S^r-\widetilde\Delta,\qquad
\Psi(x,y)={\left(\rule{0mm}{3.85mm}\nu(x+y)\,,\nu(x-y)\right)},
\end{equation}
where $\nu$ is the normalization
map at the end of the proof of Proposition~\ref{symigualdad}, 
sends relations~(\ref{reltcs})
into relations~(\ref{relemb}) and vice versa. Moreover, the resulting maps
$\Psi'\colon E_2\to \P^r\times\P^r-\Delta_{\P^r}$ and $\Psi''\colon 
\P^r\times\P^r-\Delta_{\P^r} \to E_2$ are easily seen to be equivariant with
respect to the $\mathbb{Z}/2$-action on $E_2$ 
coming from $\pi_2$ in the right triangle
of~(\ref{triangulo1}), and on $\P^r\times\P^r-\Delta_{\P^r}$ coming from
interchanging coordinates.

\medskip The result is now a direct consequence of~(\ref{class}).
\end{proof}

\begin{proof}[Proof of Theorem~\ref{Tfinal}]
Use, in this order, Definition~\ref{TCSdef}, Proposition~\ref{symigualdad},
Proposition~\ref{trick}, and~(\ref{class}).
\end{proof}

\section{Haefliger's metastable range}\label{metastable}
Most cases in Theorem~\ref{mainT} will follow directly from 
Theorem~\ref{Tfinal} 
and the following characterization of smooth 
embeddings~(proved in~\cite[Th\'eor\`eme~$1'$]{Haefligerstable}).

\begin{teorema}[Haefliger]\label{haefligerbasic}
Let $2m\ge 3(n+1)$. For a smooth compact $n$-dimensional manifold $M$,
there is a surjective map from the set of isotopy classes of smooth embeddings
$M\subset\mathbb{R}^m$ onto the set of $\,\mathbb{Z}/2$-equivariant 
homotopy classes of maps $M^*\to S^{m-1}$. Here $\mathbb{Z}/2$ acts antipodally
on $S^{m-1}$, and by interchanging coordinates on $M^*=M\times M-\Delta_M$.
\hfill\cajita
\end{teorema}

Of course, all we need from Theorem~\ref{haefligerbasic} is the fact that, 
under the stated hypothesis (the so-called metastable range), the existence 
of a smooth embedding 
$M\subset\mathbb{R}^m$ is equivalent to the existence of a 
$\mathbb{Z}/2$-equivariant map $M^*\to S^{m-1}$. 
Although not relevant for our \rojo{immediate}
purposes, it is worth remarking that the 
surjective map in Theorem~\ref{haefligerbasic} 
is explicit (see~(\ref{haefigersmap})), 
and that it is in fact bijective when $2m> 3(n+1)$.

\begin{proof}[Proof of Theorem~\ref{mainT}]
An immediate consequence of Proposition~\ref{improved} below
is that the cases with $r\ge8$ in
Theorem~\ref{mainT} lie within the metastable range hypothesis in 
Theorem~\ref{haefligerbasic}. Therefore, in those cases, 
Theorem~\ref{mainT} follows from Theorems~\ref{Tfinal}
and~\ref{haefligerbasic}. The few cases in Theorem~\ref{mainT} outside  
Haefliger's metastable range ($r\in\{1,2,4\}$) will be handled at 
the end of this section. (It would be  
interesting to know whether any one of the remaining cases 
$r\in\{3,5,6,7,10,11,12,14,15\}$ gives an actual exception to 
Theorem~\ref{mainT}.)
\end{proof}

Recall that $\P^r\times \P^r\to \P^s$ is said to be an axial map if it
is homotopically non-trivial over each axis.
From~\cite[Lemma~2.1]{AGJ:axial} we know that, when $r>15$, 
an axial map $\P^r\times 
\P^r\rightarrow \P^s$ can exist only for $2s>3r$ (in view of 
Theorem~\ref{converse} below, such an axial map that in addition was 
{\it symmetric} would yield an embedding within Haefliger's 
metastable range). We will need to consider the following slight improvement.

\begin{proposicion}\label{improved}
{For $r\in\{8,9,13\}$ or $r>15$, 
an axial map $\P^r\times\P^r\to\P^s$ can exist only
when $2s\ge3(r+1)$.}
\end{proposicion}

\begin{proof}
The main result in~\cite{AGJ:axial} (see Theorem~\ref{AGJpaper} below)
implies that the axial map hypothesis can be replaced 
by an immersion $\P^r\looparrowright\mathbb{R}^s$, and we need to prove
that, for $r$ as stated, the smallest such $s$ satisfies $2s>3r+2$.
Cases with $r\in\{8,9,13\}$ follow from inspection of~\cite{tablas}.
For $r>15$ we revisit the argument in the proof of~\cite[Lemma~2.1]{AGJ:axial}.
Pick $\rho\ge 4$ with $2^\rho\le r <2^{\rho+1}$. Each of the cases
\begin{itemize}
\item $r\le 2^\rho+3$
\item $r=2^{\rho+1}-1$
\item $2^\rho+2^{\rho-1}+2\le r \le 2^{\rho+1}-3$
\end{itemize}
can be dealt with by the corresponding non-immersion result stated 
in~\cite{AGJ:axial}.

Assume $2^\rho+4\le r \le 2^\rho+2^{\rho-1}+1$ and choose 
$\sigma\in\{1,2,\ldots,\rho-2\}$ with $2^\rho+2^\sigma+2\le r\le 
2^\rho+2^{\sigma+1}+1$. From~\cite{Davisstrongimmersion},
$\P^{2^\rho+2^\sigma+2}$ does not immerse in Euclidean space
of dimension $2^{\rho+1}+2^{\sigma+1}-4$. Therefore, in the optimal
immersion $\P^r\looparrowright\mathbb{R}^s$, we must have
$s\ge 2^{\rho+1}+2^{\sigma+1}-3$, and this easily yields the required
inequality $2s>3r+2$ when $\sigma\ge 3$ or $\rho\ge 5$. For the smaller 
cases with $\rho=4$ and $1\le\sigma \le 2$, the required
$2s>3r+2$ follows, as above, from direct inspection of~\cite{tablas}.

It remains to consider the case $r=2^{\rho+1}-2$. As the case 
$\rho=4$ follows again from inspection of~\cite{tablas}, we assume further 
$\rho\ge 5$. Let $m=2^{\rho-1}+2^{\rho-2}+2^{\rho-3}$. 
From~\cite{Davisstrongimmersion} we know that $P^{2(m+\alpha(m)-1)}$
does not immerse $\mathbb{R}^{4m-2\alpha(m)}$, where $\alpha(m)$
is the number of ones appearing 
in the dyadic expansion of $m$.  Therefore, in the optimal
immersion $\P^r\looparrowright\mathbb{R}^s$, we must have
$s\ge 4m-2\alpha(m)+1=2^{\rho+1}+2^\rho+2^{\rho-1}-5$, from which
one easily deduces the required inequality $2s>3r+2$.
\end{proof}

We close this section by describing 
%the best upper and lower bounds we 
what we know about the numerical value of $\mbox{TC}^S(\P^r)$
for the small values of $r$ not covered by the hypothesis of 
Proposition~\ref{improved}, i.e., when Haefliger's metastable range hypothesis
in Theorem~\ref{mainT} might fail to hold.

\medskip
The starting point is
\begin{equation} \label{chaineq}  
\mbox{TC}(\P^r)\le \mbox{TC}^S(\P^r)\le E(\P^r)+1.
\end{equation}
The first inequality has been proved (for any space, and not only for 
projective spaces) in~\cite[Corollary~9]{FGsymm}, whereas the second 
inequality holds without restriction on $r$ in view of 
Theorem~\ref{Tfinal} and since the construction of the map in 
Theorem~\ref{haefligerbasic} makes no use of Haefliger's range 
(see~(\ref{haefigersmap}) below). In our current range $r\le 15$, the 
precise numeric value of the 
lower bound for $\mbox{TC}^S(\P^r)$ coming from the first inequality 
in~(\ref{chaineq}) is determined by~\cite{tablas,FYT}. 
As for the upper bound, note that
the term $E(\P^r)$ on the right hand side of~(\ref{chaineq})
can even be replaced by the potentially smaller 
$E_{\hspace{.25mm}\text{TOP}}(\P^r)$,
the dimension of the smallest Euclidean space where $\P^r$ admits a 
{\em topological} embedding. Indeed, such an embedding $g\colon \P^r
\hookrightarrow\mathbb{R}^d$ determines a $\mathbb{Z}/2$-equivariant map
$\widetilde{g}\colon\P^r\times\P^r-\Delta_{\P^r}\to S^{d-1}$ by the usual 
formula
\begin{equation}\label{haefigersmap}
\rojo{\widetilde{g}(a,b)=\frac{g(a)-g(b)}{\left|g(a)-g(b)\right|}.}
\end{equation}

Now, the low dimensional cases under consideration either have
$r\le 7$ or $r\in\{10,11,12,14,15\}$. For $r=7$ and $r=15$ the use of 
$E_{\hspace{.25mm}\text{TOP}}(\P^r)$ gives more 
accurate information than that available for 
$E(\P^r)$. Indeed, the PL embeddings $\P^7\hookrightarrow\mathbb{R}^{10}$
and $\P^{15}\hookrightarrow\mathbb{R}^{23}$ constructed in~\cite{rees}
improve by 2 and 1 units, respectively, the upper bound in~(\ref{chaineq}) 
obtained from the best smooth embedding results currently known 
(see~\cite{tablas}). Also worth noticing is the fact that Rees' upper bound
$$\mbox{level}(\P^6\times\P^6-\Delta_{\P^6},\mathbb{Z}/2)\le 9,$$
obtained in~\cite[Corollary~11]{rees2},
improves by 2 units the upper bound in~(\ref{chaineq}).
On the other hand, it is elementary to check that, in all cases with 
$r\le 7$,~\cite[Theorem~17]{FGsymm} improves by one unit the 
lower bound in~(\ref{chaineq}).

\medskip
Table~1 summarizes the resulting improved bounds $\ell(r)\le \mbox{TC}^S
(\P^r)\le u(r)$. Note that the case $r=4$ does lie within Haefliger's 
metastable range, and together with the cases $r=1,2$ gives the three 
missing instances in the proof of Theorem~\ref{mainT}.

\begin{table}
\centerline{
\begin{tabular}{|c|c|c|c|c|c|c|c|c|c|c|c|c|}\hline
$r$ & $1$ & $2$ & $3$ & $4$ & $5$ & $6$ & $7$ & $10$ & $11$ 
& $12$ & $14$ & $15$ \\ \hline
$u(r)$ & $3$ & $5$ & $6$ & $9$ & $10$ & $10$ & $11$ & $18$ & 
$19$ & $22$ & $24$ & $24$ \\ \hline
$\ell(r)$ & $3$ & $5$ & $5$ & $9$ & $9$ & $9$ & $9$ & $17$ & 
$17$ & $19$ & $23$ & $23$ \\ \hline
\end{tabular}}
\caption{Upper and lower bounds for $\mbox{TC}^S(\P^r)$ for low values of $r$}
\end{table}

\begin{ejemplos}\label{bounded}{\em
\colorado{Let $\delta=(0,1,2)$,  
$i\ge (1,3,4)$, and $r=2^i+\delta$. According to~\cite{tablas}
and Theorems~\ref{mainT} and~\ref{recuperado} (below), we have
\begin{equation}\label{diferencia}
\mbox{TC}^S(\P^r)-\mbox{TC}(\P^r)=(1,2,1).
\end{equation}
This situation contrasts with the fact, proved in Section~\ref{cx}, 
that the first inequality in~(\ref{chaineq}) 
becomes an equality for all complex projective spaces.
Actually, in view of the main result in~\cite{HS}, 
there is even a (weird) possibility that the left hand side 
of~(\ref{diferencia}) might turn out
to be an unbounded function of $r$ (compare to~\cite[Example~28]{FGsymm}).
Other situations with a 
behavior resembling that in~(\ref{diferencia}) are given 
by spheres: according to~\cite{F1,FGsymm},
$\mbox{TC}^S(S^r)-\mbox{TC}(S^r)=0$ for $n$ even, 
whereas $\mbox{TC}^S(S^r)-\mbox{TC}(S^r)=1$ for $n$ odd.}
}\end{ejemplos}

\section{Relation to axial maps}\label{relax}
In Sections~\ref{mainsec} and~\ref{mainproof} 
we made it clear that, in characterizing the embedding dimension 
for real projective spaces, one might prefer to avoid the use of 
symmetric axial maps. However, in this section we analyze the way 
Theorem~\ref{mainT} is related to such maps.

\medskip
Our justification for including this section is three-fold. First, it
shows how our proof of Theorem~\ref{mainT} arose
(compare the map $\Psi$ in~(\ref{Figrande}) below with that 
in~(\ref{extendida})). Second, it illustrates the use (and, as observed 
in the next section, the limitations) of the constructions 
in Subsection~\ref{ss1.3} when applied to the case of even-torsion lens spaces.
And third, this material will allow us to make explicit comparisons
with the maps arising in the next section
(e.g., (\ref{relationse}),~(\ref{embedding}), and~(\ref{tce}) 
as generalized forms of~(\ref{relationse1}))
towards a characterization of the symmetric topological complexity of lens 
spaces (Theorem~\ref{caracte} and Proposition~\ref{symigualdade}).

\medskip
The section has been divided into three short subsections.
We start with a brief review of the axial map 
interpretation for the immersion problem of real projective spaces
(Subsection~\ref{ss1.1}), and the corresponding (partial)
interpretation known before this paper for embeddings (Subsection~\ref{ss1.3}). 
The main goal then is to compare our methods in Section~\ref{mainproof}
to those in~\cite{BFG} (Subsection~\ref{ss1.3}) 
and~\cite{FYT} (Subsection~\ref{ss1.4}).

\medskip 
There are no new results in this 
section; instead, it has a retrospective
flavor, written much in the way the ideas in this paper originally arose. The 
reader interested in our analysis and results on the symmetric topological 
complexity of lens spaces and complex projective spaces can safely skip this 
section, and proceed directly to the final 
Sections~\ref{lens} and~\ref{cx}, respectively.

\subsection{Axial maps, immersions, and topological complexity}
\label{ss1.1}
Activities were launched with Hopf's early work \cite{H} 
constructing, for $n > r $, 
a Euclidean $n$-dimensional embedding for $\P^r$ from a given symmetric 
nonsingular bilinear map
$\alpha\colon\mathbb{R}^{r+1}\times\mathbb{R}^{r+1}\to\mathbb{R}^{n+1}$. 
By restricting to unit vectors (and normalizing), this yields a symmetric
$\mathbb{Z}/2$-biequivariant 
map $\widetilde{\alpha}\colon S^r\times S^r\to S^n$,
i.e., one satisfying conditions~(\ref{relationse1}) below. Note that
$\widetilde{\alpha}$ covers an axial map $\,\widehat{\alpha}\colon \P^r\times 
\P^r\to \P^n$ that, in addition, is symmetric in the sense that the relation
$\widehat{\alpha}(a,b)=\widehat{\alpha}(a,b)$ holds for $a,b\in \P^r$. 

\medskip Using Hirsch's characterization of smooth Euclidean
immersions in terms of the geometric dimension of the normal bundle,
the relevance of (not necessarily symmetric) axial {maps} was settled
in~\cite{AGJ:axial} by showing{: 
\begin{teorema}\label{AGJpaper}
For $n>r$, the existence of 
an axial map $\P^r\times \P^r\to \P^n$ is equivalent to 
the existence of a smooth immersion $\P^r\looparrowright \mathbb{R}^n$.
{\hfill\cajita}
\end{teorema}
}

The hypothesis $n>r$ is needed only for $r=1,3,7$. 
In those cases $\P^r$ is parallelizable and  
has an optimal Euclidean immersion in codimension 1; 
however the complex, quaternion, and octonion 
multiplications yield axial maps with $n=r$.

\medskip
But the connection with robotics was established
after 30 years with M.~Farber's work (initiated in~\cite{F2,F1}) 
on the motion planning problem. The main result in~\cite{FYT} is:

\begin{teorema}\label{recuperado}
%Namely, as shown in~\cite{FYT}, 
For $r\neq 1,3,7$, TC$(\P^r)$ is the smallest
integer $n$ such that there is an axial map $\P^r\times \P^r\to \P^{n-1}$.
Consequently, TC$(\P^r)-1$ is the smallest dimension of 
Euclidean spaces where $\P^r$ can be smoothly immersed. This assertion
holds for the three exceptional values of $r$ provided 
TC$(\P^r)-1$ is replaced by TC$(\P^r)$. \hfill\cajita
\end{teorema}

\subsection{Symmetric axial maps and embeddings}\label{ss1.3}
As shown in~\cite{BFG}, the embedding problem for \rojo{real} projective spaces
can be closely modeled by keeping Hopf's original 
symmetry condition for axial maps. We give below a quick review 
of some of the main ideas in~\cite{BFG}. 

\medskip
Start by observing that a symmetric axial map 
$\widehat{\alpha}\colon \P^r\times \P^r\to \P^s$ is covered
by a map $\widetilde{\alpha}\colon S^{r}\times S^r\rightarrow S^s$ satisfying
\begin{equation}\label{relationse1}
-\widetilde{\alpha}(x,y)=\widetilde{\alpha}(-x,y)=\widetilde{\alpha}
(x,-y)\mbox{ \ and \ }\widetilde{\alpha}(x,y)=\widetilde{\alpha}(y,x)
\end{equation} 
for $x,y\in S^r$.
Under these conditions it is elementary to check 
that the composite 
\begin{equation}\label{Figrande}
V_{r+1,2}\stackrel{\Psi}{\longrightarrow}
S^{r}\times S^r\stackrel{\widetilde{\alpha}}{\longrightarrow} S^s,\qquad
\Psi(x,y)=\left(\frac{x+y}{\sqrt{2}}\,,\frac{x-y}{\sqrt{2}}\right)
\end{equation}
is a $D_4$-equivariant map. Here
$D_4$ is the dihedral group written as the wreath product
%$\mathbb{Z}/2\wr\mathbb{Z}/2$
$(\mathbb{Z}/2\times\mathbb{Z}/2)\rtimes\mathbb{Z}/2$
where $\mathbb{Z}/2$ acts on $\mathbb{Z}/2\times\mathbb{Z}/2$
by interchanging factors.
This group acts on $S^s$ via the canonical projection $(\mathbb{Z}/2
\times\mathbb{Z}/2)\rtimes\mathbb{Z}/2\to\mathbb{Z}/2$, and on
$V_{r+1,2}$ (the Stiefel manifold of orthonormal 2-frames
in $\mathbb{R}^{r+1}$)
via the restricted \rojo{left} $D_4$-action 
in $S^r\times S^r$, where $\mathbb{Z}/2\times\mathbb{Z}/2$ and 
$\mathbb{Z}/2$ act on $S^r\times S^r$ by the product antipodal-action 
and by switching coordinates, respectively.

\medskip 
{On the other hand, with the notation  
$\widetilde\Delta=\{(x,y)\in S^r\times S^r \,|\, x\neq\pm y\}$ in 
the proof of Proposition~\ref{symigualdad}, the map
$H\colon (S^r\times S^r-\widetilde\Delta)\times [0,1]
\to S^r\times S^r-\widetilde\Delta$ defined by $H(u_1,u_2,t)=
(\widetilde u_1,\widetilde u_2)$ where
\begin{eqnarray*}
\widetilde u_1 = \frac{u_1+t(v_1-u_1)}{||\, u_1+t(v_1-u_1)\,||} 
& \qquad & \hspace{-.03cm}
\widetilde u_2 = \frac{u_2+t(v_2-u_2)}{||\, u_2+t(v_2-u_2)\,||} 
\rule{0mm}{5mm}\\
v_1=w_1+w_2\hspace{1.925cm} & \qquad & v_2=w_1-w_2 \rule{0mm}{4mm}\\
w_1=\frac{u_1+u_2}{\sqrt{1+\langle u_1,u_2\rangle}} \hspace{.82cm}& \qquad & 
\hspace{-.09cm}w_2=\frac{u_1-u_2}{\sqrt{1-\langle u_1,u_2\rangle}}
\rule{0mm}{6mm}
\end{eqnarray*}
gives a $D_4$-equivariant deformation retraction of 
$S^r\times S^r-\widetilde\Delta$ onto $V_{r+1,2}$. (Figure~1 
depicts the case in which the angle between $u_1$ and $u_2$ 
is less than $90$ degrees; the situation for an angle between $90$ 
and $180$ degrees is similar, but lowering the angle to be $90$ degrees.)}
{Then, composing the retraction $H(-,1)$ with 
$\widetilde{\alpha}\circ\Psi$ and passing to 
($\mathbb{Z}/2\times\mathbb{Z}/2$)-orbit 
spaces, we get a $\mathbb{Z}/2$-equivariant map $(\P^r)^*\to S^s$.}

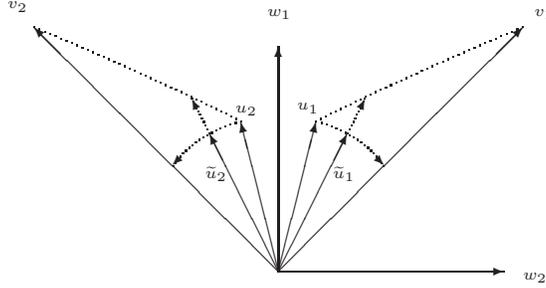
\begin{figure}
\centerline{
\begin{picture}(0,70)(0,0)
\put(0,0){\vector(0,1){60}}
\put(-3,68){\tiny$w_1$}
\put(0,0){\vector(-1,4){10}}
\put(0,0){\vector(1,4){10}}
\put(5,43){\tiny$u_1$}
\put(-11.5,42.4){\tiny$u_2$}
\put(0,0){\vector(1,0){60}}
\put(65,-2){\tiny$w_2$}
\put(0,0){\vector(1,1){65}}
\put(0,0){\vector(-1,1){65}}
\put(68,68){\tiny$v_1$}
\put(-72,70){\tiny$v_2$}
\qbezier[18](10,40)(23,36)(28,28)
\qbezier[18](-10,40)(-23,36)(-28,28)
\put(28.2,28.2){\vector(1,-1){0}}
\put(-28.2,28.2){\vector(-1,-1){0}}
\qbezier[40](10,40)(36,53)(65,65)
\qbezier[40](-10,40)(-36,53)(-65,65)
\put(0,0){\vector(1,2){18.3}}
\put(0,0){\vector(-1,2){18.3}}
\qbezier[10](18.4,36.8)(21,42)(23,46)
\qbezier[10](-18.4,36.8)(-21,42)(-23,46)
\put(23,46){\vector(1,2){0}}
\put(-23,46){\vector(-1,2){0}}
\put(14.5,25){\tiny $\widetilde{u}_1$}
\put(-19.5,25){\tiny $\widetilde{u}_2$}
\end{picture}}
\caption{The $D_4$-equivariant deformation retraction $H$}
\end{figure}

\medskip
In view of Theorem~\ref{haefligerbasic},
the above construction settles the first statement of 
Theorem~\ref{converse} below. The bulk of the work 
in~\cite{BFG} uses Haefliger and Hirsch's 
fundamental work~\cite{Haefligerstable,HH} on embeddings and 
immersions in the stable range to establish the 
second statement of Theorem~\ref{converse}.
\begin{teorema}\label{converse}
The existence of a 
symmetric axial map $\P^r\times\P^r\to\P^s$
implies the existence of a smooth embedding $\P^r\subset\mathbb{R}^{s+1}$,
provided $2s>3r$. The existence of a 
{smooth} embedding $\P^r\subset\mathbb{R}^s$
implies the existence of a symmetric axial map $\P^r\times \P^r\rightarrow 
\P^s$.\hfill \cajita
\end{teorema}
The arguments in~\cite{BFG} go a bit further.
Using the full power of 
Theorem~\ref{haefligerbasic}, it is possible to explicitly 
relate, for instance, isotopy classes of embeddings
to symmetric homotopy classes of symmetric axial maps. We
will not make use of these more complete results, though.

\medskip
Theorem~\ref{converse} can be interpreted as follows.
Let $a_S(r)$ denote the smallest
integer $k$ for which there exists a symmetric axial map
$\P^r\times \P^r\rightarrow \P^k$. It is immediate from 
Theorem~\ref{converse} that, at least for $r$ as in Proposition~\ref{improved},
\begin{equation}\label{Eas}
E(r)=a_S(r)+\delta\mbox{ \ with \ }\delta=\delta(r)\in\{0,1\}.
\end{equation}
To the best of our knowledge, the explicit value of $\delta$ (as a function 
of $r$) remains an open question. As explained in Section~\ref{mainsec},
Theorems~\ref{mainT} and~\ref{Tfinal} avoid this 
$\delta$ indeterminacy by replacing $a_S(r)$ with TC$^S(\P^r)$.

\subsection{Symmetric axial maps and symmetric TC}\label{ss1.4}
The goal of this subsection is to offer a direct comparison between 
our methods and those used in~\cite{FYT} for the non-symmetric case.
To this end, we start by observing the following obvious consequence 
of~(\ref{Eas}) and the second inequality in~(\ref{chaineq}).

\begin{corolario}\label{symmotion}
For any $r$, $\mbox{\em TC}^{S}(\P^r)\le a_S(r)+2$.\hfill\cajita
\end{corolario}

In view of Theorem~\ref{mainT}, 
this can be thought of as extending the $E$ vs.~$a_S$ relation in~(\ref{Eas})
within the topological complexity viewpoint. Of course, for $r$ as in
Theorem~\ref{mainT}, the above inequality is within one unit of being an 
equality. 

\medskip 
Corollary~\ref{symmotion} can be settled with 
a straightforward adaptation of the idea in the
first part of the proof of~\cite[Proposition~17]{FYT}. Namely, 
let $\widetilde{\alpha}
\colon S^r\times S^r\rightarrow S^s$ satisfy~(\ref{relationse1}),
with $s=a_S(r)$. 
For $1\le i\le s+1$, let $\alpha_i\colon S^r\times S^r\to \mathbb{R}$ 
be the $i$th real component of $\widetilde{\alpha}$, and set
$U_i$ to be the open subset of $\P^r\times \P^r-\Delta_{\P^r}$ consisting
of pairs $(L_1,L_2)$ of lines with 
$\alpha_i(\ell_1,\ell_2)\neq 0$, for representatives $\ell_j\in L_j\cap
S^r$, $j=1,2$. Consider the function
$s_i\colon U_i\to P_1(\P^r)$ 
defined as follows: given $(L_1,L_2)\in U_i$,
choose elements $\ell_j$ as above with $\alpha_i(\ell_1,\ell_2)>0$.
The two such possibilities $(\ell_1,\ell_2)$ and $(-\ell_1,-\ell_2)$
give the same orientation for the $2$-plane $P(L_1,L_2)$
generated by $L_1$ and $L_2$. Under these conditions,
$s_i(L_1,L_2)$ is the path rotating $L_1$ to $L_2$ in the oriented plane
$P(L_1,L_2)$.

\smallskip
Evidently  $s_i$ is a continuous section of the fibration $\ev_1$ over $U_i$. 
It is also
$\mathbb{Z}/2$-equivariant, in view of the last condition 
in~(\ref{relationse1}). Therefore, it induces a corresponding (continuous)
section $\bar{s_i}$ of the fibration $\ev_2$ over the image 
of $U_i$ under the canonical 
(open) projection $\P^r\times \P^r-\Delta_{\P^r}\to B(\P^r,2)$. But
$\P^r\times \P^r-\Delta_{\P^r}$ is covered by the
$U_i$'s, so we deduce
$\gen(\ev_2)\le s+1$. Adding 1, we get the conclusion 
in Corollary~\ref{symmotion}.{\hfill\cajita}

\medskip
Next, we elaborate on the main difference between our methods
and those in~\cite{FYT}.
Let $\xi$ be the Hopf line bundle over $\P^r$ and consider
the exterior tensor product $\xi\otimes\xi$ over $\P^r\times \P^r$.
Let $I(r)$ denote
the smallest integer $k$ such that the iterated $(k+1)$-fold Whitney 
multiple of $\xi\otimes\xi$ admits a nowhere vanishing section.
Finally, let $a(r)$ denote the smallest integer $k$ for which there is
\rojo{a (perhaps non-symmetric)} 
axial map $\P^r\times \P^r\to \P^k$. The main results in~\cite{FYT},
Corollary~5 and Proposition~17, give
\begin{equation}\label{sinsym}
I(r)+1\le TC(\P^r)\le a(r)+1,
\end{equation}
an assertion a bit sharper than its symmetric analogue 
(Theorem~\ref{mainT} and Corollary~\ref{symmotion}).
The punch line then comes from the classical fact that, for $n\neq 1,3,7$, 
both $I(r)$ and $a(r)$ agree with the dimension of the smallest Euclidean 
space where $\P^r$ admits a smooth immersion (in the case of $I(r)$ see, 
for instance, the proof of~\cite[Proposition~2.7]{TC}). However, there is no
sectioning-Whitney-multiples interpretation available 
for a symmetric version of~(\ref{sinsym}). Instead, as detailed below,
the solution comes from an adaptation of the ideas in~\cite{BFG}. 

\medskip
Recall the Borel construction $S^r\times_{\mathbb{Z}/2} S^r$ 
introduced in Section~\ref{mainproof}.
The $2$-fold  Cartesian product of the canonical projection $S^r\to \P^r$
factors through $S^r\times_{\mathbb{Z}/2} S^r$ yielding
the $\mathbb{Z}/2$-covering space
$\pi\colon S^r\times_{\mathbb{Z}/2} S^r\to \P^r\times \P^r$
in~(\ref{triangulo}).
It is well known that $\pi$ is the sphere
bundle associated to $\xi\otimes\xi\to\P^r\times \P^r$
(see for instance~\cite[Lemma~3.1]{Steer}). The relevance of such
an observation comes from~\cite[Theorem~3, and final remarks in 
Chapter~II]{schwarz}, which affirms that the Whitney multiple 
$k(\xi\otimes\xi)$ admits a global nowhere zero section for $k=\gen(\pi)$, that
is, $I(r)+1\le\gen(\pi)$.
%Now, using~\cite[Theorem~3, and final remarks in Chapter~II]{schwarz}, 
%we know that the Whitney multiple $k(\xi\otimes\xi)$ admits a global 
%nowhere zero section for $k=\gen(\pi)$. The relevance of these 
%observations comes from the standard fact that the smallest such $k$
%agrees with the smallest $n$ for which there is an axial map
%$\P^r\times \P^r\to \P^{n-1}$. 
In these terms, the work in~\cite{FYT} for settling the first inequality 
in~(\ref{sinsym}) comes from observing that the topological complexity of 
$\P^r$ is bounded from below by $\gen(\pi)$. (The second inequality 
in~(\ref{sinsym}) is actually settled in~\cite{FYT} with a
sharpening of the argument we gave above for proving
Corollary~\ref{symmotion}.) This lower bound is easily settled
in~\cite[Theorem~3]{FYT} from diagram~(\ref{triangulo}).
In fact, as a byproduct of the methods in~\cite{FYT}, it is known that
\begin{equation}\label{igualdad}
\mbox{TC}(\P^r)=\gen(\pi).
\end{equation}
But, in order to get up to this key stage in the symmetric situation, 
we needed to adapt the main trick in~\cite{BFG}. Indeed, as shown in 
Section~\ref{mainproof}, Propositions~\ref{symigualdad} and~\ref{trick}
provide us with the needed substitute for~(\ref{igualdad}), from which the 
proof of Theorems~\ref{mainT} and~\ref{Tfinal} easily follows.

\section{Lens spaces}\label{lens}
Unlike the case of real projective spaces,
the symmetric topological complexity of a lens space 
\rojo{$L^{2n+1}(m)$---the orbit space of the standard $\mathbb{Z}/m$-action on
$S^{2n+1}$---is in general not} related to its embedding dimension 
\rojo{nor, for that matter, to the level of the 
switching involution on $L^{2n+1}(m)
\times L^{2n+1}(m)-\Delta_{L^{2n+1}(m)}$. And here is an extreme example:
while the high-torsion lens spaces $L^{2n+1}(m)$ in Subsection~\ref{ss3.3} 
have symmetric 
topological complexity equal to $4n+\epsilon$, 
with $\epsilon\in\{2,3\}$ (see~(\ref{challenge})),
the level of the corresponding switching involution is at most $2n+3$, 
for all odd $m$ (see~\cite{rees2}).}

\medskip
In retrospect, the problem arises \rojo{(in the 2-local situation)} 
from the fact that the (non-symmetric) topological complexity
of $L^{2n+1}(2^e)$ actually differs from the immersion 
dimension of this manifold, and
the difference gets larger as $e$ increases, until it attains a certain
stable value (see Remark~\ref{difstable}). 

\medskip
Following the non-symmetric lead, 
in this section we (a) indicate how one can 
characterize \rojo{the symmetric topological complexity
of an even-torsion lens space
(Theorem~\ref{caracte} and Proposition~\ref{symigualdade})}, and (b) 
point out concrete differences with respect 
to a similar characterization for \rojo{its}
embedding dimension. To better appreciate the picture,
it will be convenient to start by making a summary of, and comparing to, 
the known situation in the 
non-symmetric case.

\subsection{$e$-axial maps, immersions, 
and embeddings of $L^{2n+1}(2^e)$}\label{ss3.1}
The well known relation (Theorem~\ref{AGJpaper}) between 
Euclidean immersions of real projective spaces
and (not necessarily symmetric) axial maps 
has been generalized in~\cite{ADG} 
for $2^e$-torsion lens spaces
to prove that, with the possible exceptions of $n=2,3,5$, the existence 
of an immersion
$L^{2n+1}(2^e)\immto\mathbb{R}^m$ is equivalent to the 
existence of an 
$e$-axial map $\P^{2n+1}\times_{\mathbb{Z}/2^{e-1}} 
\P^{2n+1}\rightarrow
\P^{m}$, that is, a map \rojo{that yields} a
standard 
axial map when precomposed with the canonical 
projection $\P^{\rojo{2n+1}} \times \P^{\rojo{2n+1}}\rightarrow 
\P^{\rojo{2n+1}} \times_{\mathbb{Z}/2^{e-1}} \P^{\rojo{2n+1}}$. Here the 
notation $\P^{\rojo{2n+1}} \times_{\mathbb{Z}/2^{e-1}}\P^{\rojo{2n+1}}$
refers to the usual 
Borel construction with respect to
the standard free $\mathbb{Z}/2^{e-1}$-action on $\P^{2n+1}$ with orbit space
$L^{2n+1}(2^e)$.

\medskip
At the level of covering spaces, the $e$-axial map condition translates into
having a map $\widetilde\alpha\colon S^{2n+1}\times S^{2n+1}
\rightarrow S^{m}$ satisfying
the relations
\begin{equation}\label{relationse}
\widetilde\alpha(\omega x,y)=\widetilde\alpha(x,\omega y)\mbox{ \ and \ }
\widetilde\alpha(-x,y)=-\widetilde\alpha(x,y)
\end{equation}
for $x,y\in S^{2n+1}$ 
and $\omega\in\mathbb{Z}/2^e\subset S^1$---these correspond 
to the first group of
conditions in~(\ref{relationse1}). Our first objective is to indicate how
the slight variation in~(\ref{embedding}) below of the 
obvious symmetrization of these conditions
describes the Euclidean embedding dimension
for \rojo{(arbitrary-torsion)} lens spaces.

\medskip For an integer $m\ge2$,
the product action of $\mathbb{Z}/m\times \mathbb{Z}/m$ 
on the Cartesian product
$S^{2n+1}\times S^{2n+1}$ extends to a left action of the wreath
product $G_{\rojo{m}}=(\mathbb{Z}/\rojo{m}\times \mathbb{Z}/\rojo{m})\rtimes
\mathbb{Z}/2$, where $\mathbb{Z}/2$ acts on
$S^{2n+1}\times S^{2n+1}$ by interchanging axes.
This action is stable on the
orbit configuration space $F_{\mathbb{Z}/\rojo{m}}
(S^{2n+1},2)$ consisting of
pairs in $S^{2n+1}\times S^{2n+1}$ generating different 
$\mathbb{Z}/\rojo{m}$-orbits (this is the obvious generalization 
of the space $S^r\times S^r-{\widetilde\Delta}\,$
found in Subsection~\ref{ss1.3} as well as in the proofs of 
Propositions~\ref{symigualdad} and~\ref{trick}). 
The quotient $F_{n,m}=F_{\mathbb{Z}/m}
(S^{2n+1},2)/(\mathbb{Z}/m
\times \mathbb{Z}/m)$ has an involution induced
by the action of $G_{m}$ on the orbit configuration space, and this
gives a $\mathbb{Z}/2$-equivariant model for 
$L^{2n+1}(m)\times L^{2n+1}(m)-\Delta_{L^{2n+1}(m)}$, where
$\mathbb{Z}/2$ acts by switching coordinates. In these terms, 
Theorem~\ref{haefligerbasic}
%Haefliger's 
%characterization~\cite[Th\'eor\`eme~$1'$]{Haefligerstable} (in the 
%stable range) of Euclidean embeddings for $L^{2n+1}(2^e)$
translates into:

\begin{lema}\label{haefliger}
Assume $\rojo{k}\ge 3(n+1)$. $L^{2n+1}(m)$ can be smoothly embedded in 
$\mathbb{R}^{k}$ 
if and only if there is a $\mathbb{Z}/2$-equivariant map 
$F_{n,m}\rightarrow S^{k-1}.$\hfill\cajita
\end{lema}

Of course, having a $\mathbb{Z}/2$-equivariant map 
as above is equivalent to having
a $G_{\rojo{m}}$-equivariant map $\widetilde\alpha\colon 
F_{\mathbb{Z}/m}(S^{2n+1},2)\rightarrow S^{k-1}$,
where $G_{m}$ acts on $S^{k-1}$ via the canonical projection 
${(\mathbb{Z}/m\times \mathbb{Z}/m)\rtimes
\mathbb{Z}/2}\rightarrow \mathbb{Z}/2$.
Explicitly, $\widetilde\alpha$ must satisfy
\begin{equation}\label{embedding}
\widetilde\alpha(\omega x,y)=\widetilde\alpha(x,y)=
\widetilde\alpha(x,\omega y)\mbox{ \ and \ } 
\widetilde\alpha(x,y)=-\widetilde\alpha(y,x)
\end{equation}
for $x,y\in S^{2n+1}$ 
and $\omega\in\mathbb{Z}/m\subset S^1$.
As shown in Subsection~\ref{ss1.3}, in 
the case $m=2$, the key connection between~(\ref{embedding}) and 
the symmetrized version of~(\ref{relationse}) 
is given by the ideas in~\cite{BFG}, which teach us how to 
take care of the deleted ``equivariant diagonal'' 
in $F_{\mathbb{Z}/2}(S^{2n+1},2)$. 
Unfortunately, we have 
not succeeded in obtaining such a connection for larger values of $m$.
The major problem seems 
to be given by the apparent lack\footnote{\rojo{Using the gradient flow 
navigation technique in~\cite[Section~4.4]{farberinvitation}, 
Armindo Costa's current Ph.D.~work at Durham University offers a nice 
explanation of the way this problem arises.}} 
of a suitable equivariant deformation retraction of 
$L^{2n+1}(m)\times L^{2n+1}(m)-\Delta_{L^{2n+1}(m)}$
that plays the  role 
of $V_{2n+2,2}$ in the $m=2$ arguments of~\cite{BFG} 
described in Subsection~\ref{ss1.3}. It is worth remarking that,
in the symmetric $m=2$ situation of Section~\ref{mainproof}, we do 
make an indirect use---through the map $\Psi$ in~(\ref{extendida})---of 
this equivariant deformation retraction.
This problem will reappear, in a slightly different
form, in regard to a potential characterization for 
the symmetric topological complexity of $m$-torsion lens spaces
in terms of the $\mathbb{Z}/m$-biequivariant maps of the next 
subsection (see Remark~\ref{nice} {below}). 

\subsection{Symmetric biequivariant maps and 
TC$^S$ of lens spaces}\label{ss3.2}
As shown in~\cite{TC}, when $m$ is even
the (non-symmetric) topological complexity of $L^{2n+1}(m)$ 
turns out to be 
(perhaps one more than) the smallest odd integer $k$ for which there is 
a $\mathbb{Z}/\rojo{m}$-biequivariant 
map $\widetilde\alpha\colon S^{2n+1}\times 
S^{2n+1}\rightarrow S^{k}$, that is, a map satisfying
the (stronger than~(\ref{relationse})) requirements
$$
\widetilde\alpha(\omega x,y)=\widetilde\alpha(x,\omega y)=
\omega \,\widetilde\alpha(x,y),
$$
for $x,y\in S^{2n+1}$ and $\omega\in\mathbb{Z}/\rojo{m}\subset S^1$. 
Alternatively,
if $c\colon S^{2n+1}\to S^{2n+1}$ stands for complex conjugation 
in every complex coordinate, then by precomposing with $1\times c$,
a $\mathbb{Z}/m$-biequivariant map as above can equivalently 
be defined through the requirements
\begin{equation}\label{tce}
\widetilde\alpha(\omega x,y)=\omega\, \widetilde\alpha(x,y)=
\widetilde\alpha(x,\omega^{-1}y).
\end{equation}

In analogy to the $a_S$ notation introduced at the end of 
Subsection~\ref{ss1.3} to measure the existence of  symmetric axial maps,
the following definition {(which, up to composition with
$1\times c$, corresponds to the symmetrized version of the number 
\rojo{$s_{n,m}$} defined in~\cite{TC})}
is intended to measure the existence of symmetric 
$\mathbb{Z}/m$-biequivariant maps.

\begin{definicion}\label{biequivariante}{\em
For integers $n$ and \rojo{$m$}, let \rojo{$b^S_{n,m}$} denote the smallest
integer $k$ such that there is a map
$\widetilde\alpha\colon S^{2n+1}\times 
S^{2n+1}\rightarrow S^{2k-1}$ satisfying~(\ref{tce}) and
\begin{equation}\label{simetrizando}
\widetilde\alpha(x,y)=\widetilde\alpha(y,x)
\end{equation}
for $x,y\in S^{2n+1}$ and $\omega\in\mathbb{Z}/m\subset S^1$.
}\end{definicion}

The next result gives 
our characterization for (half the value of) 
the symmetric topological complexity of 
\rojo{even-torsion} lens spaces.
The proof will be postponed to the end of the subsection.

\begin{teorema}\label{caracte}
\rojo{For even $m$}, the integral part of 
$\frac{1}{2}\mbox{\em TC}^S(L^{2n+1}(m))$ agrees with
the smallest integer \rojo{$k$} such that
there is a map $\widetilde\alpha\colon F_{\mathbb{Z}/m}(S^{2n+1},2)
\to S^{2k-1}$ satisfying~{\em (\ref{tce})} and~{\em (\ref{simetrizando})}
for $x,y\in S^{2n+1}$ and $\omega\in\mathbb{Z}/m\subset S^1$.
\hfill\cajita
\end{teorema}

Most of the work in~\cite{FGsymm} goes in the direction of giving 
strong lower bounds for TC${}^S$. However, there seems to be a relative lack of 
suitable upper bounds; the only ones\footnote{\rojo{The relativized notion of 
topological complexity~\cite[Sections~4.3 and~4.4]{farberinvitation} 
seems to lead to new such upper bounds.}} 
we are aware of are derived, some way
or other, from Schwarz's general estimate for the genus of a fibration
$F\to E\to B$ in terms of the dimension of $B$ and the connectivity of $F$
(\cite[Theorems~$5$ and $5'$]{schwarz}). For instance, 
in~\cite[Proposition~10]{FGsymm} the upper bound
\begin{equation}\label{diml}
\mbox{TC}^S(M)\le 2d+1
\end{equation}
is derived for any $d$-dimensional closed smooth manifold $M$.
In the case $M=L^{2n+1}(m)$, Corollary~\ref{cota} below (which is an
immediate consequence of Theorem~\ref{caracte}) offers an alternative
to~(\ref{diml}) that takes not only dimension into account, but also torsion.
Theorem~\ref{stable1} below
gives a typical example (in the non-symmetric setting, though) 
of the potential use of this kind of result.

\begin{corolario}\label{cota}
\rojo{For even $m$,} the integral part of $\frac{1}{2}
\mbox{\em TC}^S(L^{2n+1}(m))$ 
is no greater than $b^S_{n,m}$.\hfill\cajita
\end{corolario}

\begin{nota}\label{nice}{\em 
In the direction of exploring a possible 
symmetric analogue of the main result 
in~\cite{TC}, it would be {useful} to make precise how much 
the above upper bound differs from being an equality.
The main obstruction to such a goal 
seems to be the apparent lack of an analogue for lens spaces
of the map $\Psi$ in~(\ref{extendida}) and~(\ref{Figrande}).
}\end{nota}

We close this subsection with the proof of Theorem~\ref{caracte}. 
As will quickly become clear, the details are formally the same as
in the $m=2$ case.
The $m$-analogue of~(\ref{triangulo}), 
first considered in~\cite[Theorem~3]{FYT}, reads
$$
\begin{picture}(0,39)(130,17)
\put(0,51){
\xymatrix{
{ P(L^{2n+1}(m)) } \ar[rr]^f \ar[rd]_{\ev} & 
& {\left(\rule{0mm}{3.85mm}S^{2n+1}\times S^{2n+1}\right)\left/ 
\rule{0mm}{3.87mm}\,\mathbb{Z}/m
\right.} \ar[ld]^{\pi} \\
 & { L^{2n+1}(m) \times L^{2n+1}(m) }  & 
}}\end{picture}
$$
The orbit space in the upper right corner is taken with respect to the 
diagonal $\mathbb{Z}/m$-action. The map $f$, whose definition is the 
obvious generalization
of that in the case $m=2$, is $\mathbb{Z}/2$-equivariant. 
In these conditions, the analogue of~(\ref{triangulo1}) and the proof of 
Proposition~\ref{symigualdad} generalize in a straightforward way to 
produce the following characterization of TC$^S(L^{2n+1}(m))$.

\begin{proposicion}\label{symigualdade}For even $m$,
$$\mbox{\em TC}^S\left(L^{2n+1}(m)\right)-1=
\gen\left(\rule{0mm}{3.7mm}\pi_{2,m}\colon E_{2,m}\longrightarrow
B(L^{2n+1}(m), 2)\right).$$ Here $E_{2,m}$ is the quotient of 
$F_{\mathbb{Z}/m}(S^{2n+1},2)$ by the \rojo{two} 
relations $(x,y)\sim(\omega x, 
\omega y)$ and $(x,y)\sim(y,x)$.
Moreover, $\pi_{2,m}$ is a $\mathbb{Z}/m$-cover with 
$\mathbb{Z}/m$ acting on $E_{2,m}$
as $\omega\cdot [x,y]{}={}[\omega x,y]$, for $x,y\in S^{2n+1}$ 
and $\omega\in
\mathbb{Z}/m\subset S^1$, where $[x,y]$ stands for the class
of the pair $(x,y)$.\hfill\cajita
\end{proposicion}

\begin{nota}\label{par}{\em
The requirement that $m$ be even is used in the construction of
the map $g_1$ in the proof of Proposition~\ref{symigualdad}. 
Indeed, the coordinates 
of a pair $(x_1,x_2)\in F_{\mathbb{Z}/m}(S^{2n+1},2)$ cannot be antipodal when
$m$ is even, so that~(\ref{curva}) is well defined.
}\end{nota}

Theorem~\ref{caracte} is now a direct consequence of 
Proposition~\ref{symigualdade} and the 
following $m$-analogue of~(\ref{class}) 
(proved in full generality in~\cite[Corollary 1, pg.~97]{schwarz}).

\begin{lema}\label{classe}
The canonical $\mathbb{Z}/m$-cover $S^{2n-1}\to L^{2n-1}(m)$ 
classifies
$\mathbb{Z}/m$-covers of genus at most $2n$.\hfill\cajita
\end{lema}

\subsection{(Non-symmetric) TC of high-torsion lens spaces}\label{ss3.3}
For an integer $m\ge2$ 
we say that a  lens space $L^{2n+1}(m)$ is of high torsion when 
$m$ does not divide the binomial coefficient $\binom{2n}{n}$. 
A lens space that is not of high torsion will be said to be of low torsion.
The (non-symmetric) topological complexity of a 
high-torsion lens space has recently been settled 
in~\cite{FGops}.

\begin{teorema}[\rojo{Farber-Grant}]\label{stable}
For a high-torsion lens space, 
$\mbox{\em TC}(L^{2n+1}(m))=4n+2$.\hfill\cajita
\end{teorema}

\begin{nota}\label{difstable}{\em 
This result is the analogue of the following (2-local) situation.
For a fixed $n$, the immersion dimension of $L^{2n+1}(2^e)$
is a bounded non-decreasing function of $e$ which, therefore, 
becomes stable for large $e$. As explained in~\cite{nonimm}
and~\cite[Section~6]{TC}, the stable value of the immersion dimension
is expected to be attained roughly \rojo{when 
$2^e$ does not divide $\binom{2n}{n}$}---with an expected value close to
the immersion dimension of the complex projective $n$-dimensional space.
{A very concrete situation, which compares TC to the immersion dimension
of lens spaces, is illustrated in Example~\ref{numerico} below.}
}\end{nota}

The converse implication in the statement of Theorem~\ref{stable}
is true when $m$ is even. In fact, we extend Farber-Grant's result
to the first case outside the high-torsion 
range by combining the techniques in~\cite{FGops}
with the $\mathbb{Z}/m$-biequivariant map
characterization of TC$(L^{2n+1}(m))$ 
discussed at the beginning of Subsection~\ref{ss3.2}. The result arose
from an e-mail 
exchange, dating back to mid 2007, 
between the first author and Professor Farber 
concerning the results in~\cite{FGops}.

\begin{teorema}\label{stable1}
\rojo{Let $m$ be even. If $L^{2n+1}(m)$ is of low torsion, then
$\mbox{\em TC}(L^{2n+1}(m))\le 4n$, with equality when $m$ does not divide
$\binom{2n-1}{n}$}.
\end{teorema}

\begin{proof}
Proposition~2.2 and Theorem~2.9 in~\cite{TC} yield 
$\mbox{TC}(L^{2n+1}(m))\le 4n$. The \rojo{rest comes 
from~\cite[Theorem~11]{FGops} (with $k=n$ and $\ell=n-1$)}.
\end{proof}

Since $\binom{2n}{n}=2\binom{2n-1}{n}$, the final 
part in the statement of Theorem~\ref{stable1} refers to a 2-local property
of $m$, namely, that the highest exponent of $2$ in $m$ agrees with that in
$\binom{2n}{n}$.

\begin{ejemplo}\label{2local}{\em
\rojo{It is well known that the highest power of $2$ dividing 
$\binom{2n}{n}$ is $\alpha(n)$, the number of ones in the dyadic expansion of
$n$. In particular, $\mbox{TC}(L^{2n+1}(2^e))=4n+2$ for $e>\alpha(n)$. 
Theorem~\ref{stable1} now gives 
$\mbox{TC}(L^{2n+1}(2^{\alpha(n)}))=4n$. Since 
$\mbox{TC}(L^{2n+1}(2))$ is (one more than) the immersion
dimension of the real projective space 
$L^{2n+1}(2)$~(Theorem~\ref{recuperado}),
it is highly desirable to get as much information as possible
on the value of $\mbox{TC}(L^{2n+1}(2^{e}))$ 
as $e$ goes from $\alpha(n)-1$ down to $1$.}
}\end{ejemplo}

\begin{ejemplo}\label{numerico}{\em
\rojo{Table~2} summarizes the topological complexity 
and immersion dimension for $L^{2n+1}(2^e)$ and $\mathbb{C}\P^n$
in the case $n=2^r+1$ with $r\ge 1$. The information is taken 
from~\cite{tablas,FYT} in the case of $\P^{2n+1}$, 
from~\cite{nonimm,tom} in the case of the immersion dimension of
$L^{2n+1}(2^e)$ for $e\ge 2$, from~\cite[Corollary~2]{FYT} 
in the case of TC$(\mathbb{C}\P^n)$,
and from~\cite{AH,amiya} in the case of the immersion dimension of 
$\mathbb{C}\P^n$. Note that in the case under consideration
TC$(\mathbb{C}\P^n)$ is just half the stable value of 
TC$(L^{2n+1}(2^e))$ (i.e., for $e\ge 3$). Such a behavior
comes from the fact that $\mathbb{C}\P^n$ is simply connected and from
Schwarz's estimates~\cite[Theorem~5]{schwarz} 
for the genus of a fibration.}\end{ejemplo}

\begin{table}
\centerline{
\begin{tabular}{|c|c|c|c|c|}\hline
 & $\P^{2n+1}$ & $L^{2n+1}(4)$&$\rule{0mm}{4mm}L^{2n+1}(2^e)$ 
\ $e\ge 3$ & $\mathbb{C}\P^n$\\ \hline
\rule{0mm}{4mm}TC & $4n-3$ \ ($r\ge 2$) & $4n$ & $4n+2$ & $2n+1$\\
 & \rule{0mm}{4mm}$4n-4$ \ ($r= 1$) & & &\\ \hline
\rule{0mm}{4mm}Imm & $4n-4$ & $4n-3$ & $4n-2$ & $4n-3$ \\\hline
\end{tabular}}
\caption{TC vs.~Imm for $2^e$-torsion lens spaces ($n=2^r+1$, $r\ge 1$)}
\end{table}

We close this section \rojo{by proposing} 
what we believe should be an accessible
challenge: Determine the symmetric topological complexity of 
high-torsion lens spaces. We remark that, in the high-torsion range,
the inequalities 
\begin{equation}\label{challenge}
4n+2\leq \mbox{TC}^S(L^{2n+1}(\rojo{m}))\le 4n+3
\end{equation}
follow from~(\ref{diml}), Theorem~\ref{stable}, and the analogue for 
lens spaces of the first inequality in~(\ref{chaineq}).

\section{Complex projective spaces}\label{cx}
The (non-symmetric) topological complexity of 
the $n$-dimensional complex projective space 
was computed in~\cite[Section~3]{FYT} to be TC$(\mathbb{C}\P^n)=2n+1$.
In this brief final section we show that the same value holds
in the symmetric case.

\begin{teorema}\label{TCScx}
{\em TC}$^S(\mathbb{C}\P^n)=2n+1.$
\end{teorema}

\begin{proof} In view of the analogue for complex projective spaces
of the first inequality in~(\ref{chaineq}),
we only need to show that TC$^S(\mathbb{C}\P^n)\le 2n+1$. The 
diagram of pull-back squares
$$
\xymatrix{
{ P(\mathbb{C}\P^n) } \ar[d]^{\ev} &
{ P_1(\mathbb{C}\P^n) } \ar[l] \ar[r] \ar[d]^{\ev_1} 
& { P_2(\mathbb{C}\P^n) } \ar[d]^{\ev_2} \\
{\mathbb{C}\P^n\times \mathbb{C}\P^n } &
{\mathbb{C}\P^n\times \mathbb{C}\P^n -\Delta_{\mathbb{C}\P^n}} 
\ar[l] \ar[r] & {B(\mathbb{C}\P^n,2)}
}$$
where horizontal maps on the 
left are inclusions, and horizontal maps on the right
are canonical projections onto $\mathbb{Z}/2$-orbit spaces, shows that
the common fiber for the three vertical maps is the {path} connected 
space $\Omega \mathbb{C}\P^n$. In particular, Theorem~$5'$ in~\cite{schwarz} 
applied to $\ev_2$ gives
$$
\mbox{TC}^S(\mathbb{C}\P^n)=\gen(\ev_2)+1\le\frac{\dim
\left(Y\right)}{2}+2
$$
where $Y$ is any CW complex having the homotopy type of 
$B(\mathbb{C}\P^n,2)$.
The required inequality follows since, {as indicated} below, 
there is such a model $Y$ having $\dim(Y)=4n-2$.
\end{proof}
%\begin{proposicion}\label{CWtype}
%For a smooth closed manifold $M$, the unordered configuration space
%$B(M,2)$ has the homotopy type of a compact smooth manifold with boundary.
%\end{proposicion}
%
%\noindent {\it Proof. }Choose
%an open tubular neighborhood $U$ of the diagonal in $M\times M$ which is
%invariant under the $\mathbb{Z}/2$ action on $M\times M$ that switches 
%factors. This involution is closed and free on
%the smooth compact manifold (with boundary) $X=M\times M-U$. Moreover,
%$X$ is a smooth deformation retract of $M\times M-\Delta_M$; 
%indeed the deformation retraction may be taken to be equivariant with 
%respect to the involution.  Passing to orbit spaces shows that $B(M,2)$ 
%contains, as a deformation retract, the $\mathbb{Z}/2$ orbit space 
%$Y$ obtained from $X$ under the involution.  The desired conclusion follows 
%since $Y$ is also a compact smooth manifold with boundary.
%\hfill\cajita
%
%\begin{nota}\label{Yasui}{\em
%Of course, 
%any compact smooth $m$-dimensional manifold has the homotopy type of an
%$m$-dimensional CW complex (e.g., by Morse theory). 
%In the case of the smooth homotopy type
%of $B(\mathbb{C}\P^n,2)$ in Proposition~\ref{CWtype},
%an explicit model 

In the proof of~\cite[Proposition~10]{FGsymm} it is observed that,
for a smooth closed $m$-dimensional manifold $M$, $B(M,2)$ has the
homotopy type of a ($2m-1$)-dimensional CW complex. Although this is 
certainly enough for completing the proof of Theorem~\ref{TCScx},
we point out that an explicit (and smaller) model for $B(\mathbb{C}\P^n,2)$
was described by Yasui 
in~\cite[Proposition~1.6]{yasui}. We recall the details. 
The unitary group $\mathrm{U}(2)$ has the two subgroups 
\begin{itemize}
\item[$T^2$:] diagonal matrices, and 
%\vspace{-3.7mm}
\item[$G\;$:] matrices in $T^2$ together with those of the form 
$\left(\begin{array}{cc}
               0 & z_1 \\
               z_2 & 0 
          \end{array}\right)$ for $z_1,z_2\in S^1$.
%\vspace{-2mm}
\end{itemize}
Consider the standard action of $\mathrm{U}(2)$ on 
the complex Stiefel manifold $W_{n+1,2}$ of orthonormal 2-frames 
in $\mathbb{C}^{n+1}$ with quotient the Grassmann manifold
of complex 2-planes in $\mathbb{C}^{n+1}$. Yasui's model for $B(\mathbb{C}
\P^n,2)$ is the corresponding quotient $W_{n+1,2}/G$.
Note that $\dim(G)=\dim(T^2)=2$, so that the dimension of Yasui's model is
$$\dim(W_{n+1,2})-2=4n-2.$$
%}\end{nota}

\bigskip\bigskip

Jes\'us Gonz\'alez\quad {\tt jesus@math.cinvestav.mx}

%\medskip
{\sl Departamento de Matem\'aticas, CINVESTAV--IPN 

Apartado Postal 14-740 M\'exico City 07000}

\bigskip\medskip

Peter Landweber\quad {\tt landwebe@math.rutgers.edu}

%\medskip
{\sl Department of Mathematics, Rutgers University

110 Frelinghuysen Rd,  Piscataway, NJ 08854-8019}

\end{document}